\documentclass[10pt]{article}
\usepackage{amssymb}
\usepackage{amsthm}
\usepackage[utf8]{inputenc}

\def\Bbb{\mathbb}
\def\R{{\Bbb R}}
\def\Q{{\Bbb Q}}
\def\K{{ K}}

\newtheorem{lemma}{Lemma}
\newtheorem{theorem}{Theorem}
\newtheorem{cor}{Corollary}
\newtheorem{prop}{Proposition}

\begin{document}

\title{Unit distance graphs with ambiguous chromatic number}
\author{Michael S. Payne\footnote{This work began as part of a Bachelor's thesis at Monash University, and was extended while the author was studying with the support of the Berlin Mathematical School.} \\
\small Institut für Mathematik, TU Berlin\\
\small and \\
\small School of Mathematical Sciences, Monash University. \\
\small \texttt{michaelstuartpayne@gmail.com}}

\maketitle

\begin{abstract}

First L{\'a}szl{\'o} Sz{\'e}kely and more recently Saharon Shelah and Alexander Soifer have presented examples of infinite graphs whose chromatic numbers depend on the axioms chosen for set theory. The existence of such graphs may be relevant to the Chromatic Number of the Plane problem. In this paper we construct a new class of graphs with ambiguous chromatic number. They are unit distance graphs with vertex set $\R^n$, and hence may be seen as further evidence that the chromatic number of the plane might depend on set theory.

\end{abstract}

\section{Introduction}

The Chromatic Number of the Plane problem asks how many colours are required to colour the Euclidean plane if points that are distance 1 apart must receive different colours. The number is known to be between 4 and 7 inclusive. For a comprehensive history see \cite{soifbook}. We may view the problem as that of colouring an infinite graph lying in the plane. This graph, which by abuse of notation we denote $\R^2$, has all points of the plane as its vertices and edges between points that are distance 1 apart. Any graph in the plane with straight unit length edges is therefore a subgraph of $\R^2$.

In 1984 L{\'a}szl{\'o} Sz{\'e}kely investigated the difference between the usual chromatic number ($\chi$) and the measurable chromatic number ($\chi_m$) for geometric graphs, the latter being the chromatic number when only Lebesgue measurable colour sets are allowed \cite{sakey}. He gave an example of a graph which could be 2-coloured in general, but which needed 3 colours in the measurable case. It consisted of the points on the unit circle with two points joined by an edge if the arc length between them was some fixed irrational multiple of $\pi$. Sz{\'e}kely concluded that (assuming the Axiom of Choice) chromatic number and measurable chromatic number were not in general the same.

In a recent series of papers Saharon Shelah and Alexander Soifer presented some more graphs with $\chi \neq \chi_m$ \cite{ss03,ss04,soifer05}. They made the dependence on set theory more explicit by considering two systems of axioms in particular. Firstly, under the system consisting of the Zermelo-Fraenkel axioms along with the full Axiom of Choice the graphs were found to have a finite chromatic number. In each case a colouring was given that relied on the Axiom of Choice. The second system of axioms limited the Axiom of Choice to a weaker form, the Principle of Dependent Choices, and introduced an Axiom of Lebesgue Measurability. This axiom states that every subset of the real numbers is Lebesgue measurable. Under this new system the chromatic numbers of the graphs were found to be uncountable.

The two different viewpoints, one contrasting normal with measurable chromatic number, and the other considering chromatic number under two different systems of axioms, are essentially equivalent for our present purposes. Here we will follow the terminology of Sz{\'e}kely and use $\chi$ and $\chi_m$ to distinguish between the two situations. We say that a graph has \emph{ambiguous} chromatic number when $\chi \neq \chi_m$. Unless otherwise indicated, in what follows all references to `measure' and `measurability' refer to $n$-dimensional Lebesgue measure which we will denote by $\mu$.

The purpose of this paper is to present a new family of graphs with ambiguous chromatic number. Unlike Sz{\'e}kely's example, the new examples have all of $\R^n$ as their vertex set, and unlike Shelah and Soifer's graphs, they have unit length edges and finite chromatic number in both situations.

\section{The construction}

Throughout the following $\K$ will always be a field with $\Q \subseteq \K \subseteq \R$. The Euclidean metric on $\K^n$ induces a unit distance graph which we again denote $K^n$  (we also suppose $n \geq 2$ throughout).  Now we construct the graph $T_{\K^n}$ by translating the graph $\K^n$ to all points in $\R^n$. Hence the vertex set becomes $\R^n$ and two vertices are joined by an edge if their difference is a unit vector in $\K^n$. There are several values of $\K$ and $n$ for which $\chi(\K^n)$ is known. These will become important later when we discuss $\chi({T_{\K^n}})$, but first let us consider the case of measurable colourings.

\subsection{Measurable Colourings}

We begin with a few measure theoretic definitions. For a point $x \in \R^n$ and a measurable set $S \subset \R^n$ the \emph{Lebesgue density} of $S$ at $x$ is
\[ d_S(x) := \lim_{\epsilon \to 0} \frac{\mu(B_{\epsilon}(x) \cap S)}{\mu(B_{\epsilon})}.\]
We define the \emph{essential part} $\tilde{S}$ of $S$ to be the set of points where $S$ has Lebesgue density 1.

The graph $T_{\K^n}$ has two important properties. Firstly, since rational points are dense on the unit $n$-sphere (see for example \cite{schmutz}), each vertex is connected to a dense set on the unit sphere around it. Secondly, the edge set of $T_{\K^n}$ is invariant under real translations, that is, if there is an edge at one point then parallel copies exist at all other points. These properties allow us to prove the following useful lemma.

\begin{lemma}\label{useful}
Let $S$ be a measurable set which is admissible as a colour set for $T_{\K^n}$ and suppose that $x\in \R^n$ is at unit distance from a point in $\tilde{S}$. Then $d_S(x) =0$.
\end{lemma}

\begin{proof}
Take any $\delta > 0$ and suppose $x$ is at unit distance from $y \in \tilde{S}$. Then since $d_S(y)=1$ we can find $\epsilon > 0$ small enough that 
\[ \frac{\mu(B_{\epsilon}(y) \cap S)}{\mu(B_{\epsilon})} \geq 1-\delta.\]
The density of the neighbours of $y$ allows us to find a neighbour $x'$ so close to $x$ that
\[ \frac{\mu(B_{\epsilon}(x) \setminus B_{\epsilon}(x'))}{\mu(B_{\epsilon})} \leq \delta.\]
By considering translations of the edge $(x',y)$ within these neighbourhoods it is clear that
\[\frac{\mu(B_{\epsilon}(y) \cap S)}{\mu(B_{\epsilon})} + \frac{\mu(B_{\epsilon}(x') \cap S)}{\mu(B_{\epsilon})} \leq 1.\]
Combining these inequalities gives us
\[ \frac{\mu(B_{\epsilon}(x) \cap S)}{\mu(B_{\epsilon})} \leq \frac{\mu(B_{\epsilon}(x) \setminus B_{\epsilon}(x'))}{\mu(B_{\epsilon})} + \frac{\mu(B_{\epsilon}(x') \cap S)}{\mu(B_{\epsilon})} \leq 2\delta .\]
Since $\delta$ can be arbitrarily small the conclusion follows.
\end{proof}

In 1981 Falconer showed that $\chi_m(\R^n) \geq n+3$ \cite{falconer}. Our aim is to adapt his proof to show that the same holds for $T_{\K^n}$.
We will use the following two lemmas of Falconer without modification. The first was proved by Croft in \cite{croftlat}.

\begin{lemma}\label{pos} Let $B$ be a non-empty subset of $\R^n$ with $\mu(B)=0$ and $C$ be a countable configuration of points in $\R^n$. Then given a point $x \in C$ there exists a rigid motion $m$ such that $m(C) \cap B = \{ m(x) \}$. Furthermore, almost all rotations (in the sense of rotational measure) of $m(C)$ about $m(x)$ have this property.
\end{lemma}

\begin{lemma}\label{boundary} Let $S$ be a Lebesgue measurable subset of $\R^n$ with $\mu(S)>0$ and $\mu(\R^n \backslash S)>0$, then $\partial S$ is non-empty, $\mu(\partial S)=0$ and $\tilde{S}$ is a Borel set.
\end{lemma}

\noindent Hence we see that if we have a covering of the plane by measurable sets $S_1, \dots, S_k$, then $\R^n \setminus \bigcup \tilde{S}_i = \bigcup \partial S_i$ and hence has measure 0. We now have enough to prove the following.

\begin{prop}\label{usecor} The measurable chromatic number of $T_{\K^n}$ is at least the (general) chromatic number of $\R^n$. That is, $\chi_m(T_{\K^n}) \geq \chi(\R^n)$.
\end{prop}

\begin{proof} A theorem of Erd\H{o}s and de Bruijn \cite{EdB} says that $\chi(\R^n)$ is realised on a finite unit distance graph, call it $G$. Suppose we have a measurable colouring of $T_{\K^n}$ by the sets $S_1, \dots, S_k $ with $k < \chi(\R^n)$. By Lemma \ref{pos} we can place $G$ so that its vertices lie in the union of the $\tilde{S}_i$. Since $k < \chi(G)$ there must be an edge of $G$ that has both vertices in $\tilde{S}_j$ for some $j$. This is a contradiction by Lemma \ref{useful}.
\end{proof}

Finally we need a slight modification of Falconer's fourth lemma and its corollary.

\begin{lemma}\label{circ} Let $\Sigma$ be a circle in $\R^2$ of radius $r>1/2$ such that $\theta = 2 \arcsin \!\left( \frac{1}{2r} \right)$ is an irrational multiple of $\pi$. Suppose almost all the points on $\Sigma$ (in the sense of circular measure) lie in $\tilde{S}_1$ or $\tilde{S}_2$. Then at least one of $\tilde{S}_1$ or $\tilde{S}_2$ realises distance 1.
\end{lemma}

\begin{proof} The only difference is that in the conclusion $\tilde{S}_1$ or $\tilde{S}_2$ realises distance 1 instead of $S_1$ or $S_2$. This new conclusion is actually an intermediate step in Falconer's proof \cite{falconer}.
\end{proof}

\begin{cor}\label{nsphere} Let $\Sigma$ be an $(n-1)$-sphere of radius $> \frac{1}{2}$ in $\R^n$, where $n>2$. Suppose $\R^n$ is divided into measurable sets $S_i$ such that almost all points of $\Sigma$ lie in $\tilde{S}_1$ or $\tilde{S}_2$. Then at least one of $\tilde{S}_1$ or $\tilde{S}_2$ realises distance 1.
\end{cor}

\begin{proof} Take a suitable affine plane section. \end{proof}

With all this preparation we can now prove our main theorem. As always, $\Q \subseteq \K \subseteq \R$.

\begin{theorem}\label{main} Any colouring of the graph $T_{\K^n}$ by measurable sets requires at least $n + 3$ colours. That is, $\chi_m(T_{\K^n}) \geq n+3$.
\end{theorem}
\begin{proof}
Suppose we have a colouring of $T_{\K^n}$ by $n+2$ measurable sets $S_0, \dots , S_{n+1}$. As in Falconer's proof we consider a configuration $C$ of $n+2$ points $x_1, \dots , x_{n+2}$ consisting of a unit $n$-simplex formed by the points $x_1, \dots, x_{n+1}$, along with the image $x_{n+2}$ of the point $x_1$ reflected in the hyperplane containing $x_2, \dots , x_{n+1}$.
Let $B= \R^n \setminus \bigcup \tilde{S}_i$. By Lemma \ref{pos} we can place $C$ with $x_1$ in $B$, and so that for almost all rotations $\rho$ of $C$ about $x_1$ we have $\rho(C) \cap B = \{x_1\}$. We can assume that $x_1$ is in the boundary of at least two sets, say $S_0$ and $S_1$. Then for all such $\rho$ we use Lemma \ref{useful} to deduce that the $\rho(x_i)$ are in one each of the $\tilde{S}_i$ for $2 \leq i \leq n+1$, and that $\rho(x_{n+2})$ is in either $\tilde{S}_0$ or $\tilde{S}_1$. Hence we know that the $(n-1)$-sphere around $x_1$ of radius $|x_1 - x_{n+2}|$ lies almost all in $\tilde{S}_0 \cup \tilde{S}_1$. We refer to Falconer's proof for that fact that this radius satisfies the conditions of Lemma \ref{circ} in the case $n=2$, and then apply it and Corollary \ref{nsphere} and also Lemma \ref{useful} to get the result.
\end{proof}

\subsection{General Colourings}

For general colourings we have the following result.

\begin{prop}\label{ccolour} $\chi(\K^n) = \chi(T_{\K^n})$.
\end{prop}
\begin{proof}
The translates of $\K^n$ that make up $T_{\K^n}$ are disconnected from each other so each one can be coloured 
independently.
\end{proof}
\noindent The obvious way to colour each translate is by translating a fixed colouring of $\K^n$ to each one. For $\K$ countable we must apply the Axiom of Choice to an uncountable collection of sets to select representatives of the translates on which to start the colouring. If we choose the representatives from inside the unit cube then the set of representatives is a classic Vitali type non-measurable set, so the colour sets of our colouring are countable unions of non-measurable sets. It is not surprising then that such colour sets may turn out to be non-measurable.

\section{Ambiguous cases}\label{last}

Returning at last to the the topic of ambiguity, comparing Proposition \ref{ccolour} and Theorem \ref{main} we can now see that if $\chi(\K^n) < n+3$ then $T_{\K^n}$ has ambiguous chromatic number. We note that it is clear that $\chi(T_{\K^n}) \leq \chi_m(T_{\K^n}) \leq \chi_m(\R^n)$, and that $\chi_m(\R^n)$ is finite for all $n$ because the tile based colourings that establish upper bounds on $\chi(\R^n)$ are measurable colourings. So when ambiguity occurs for $T_{\K^n}$ the chromatic numbers in both cases will be finite. 

Firstly let us consider the case where $\K=\Q$. The chromatic number of $\Q^n$ has been studied quite extensively and it is known that $\chi(\Q^2)=2$, $\chi(\Q^3)=2$ and $\chi(\Q^4)=4$ (see \cite{woodall} for the first and \cite{benper} for the other two claims). Hence the chromatic number of $T_{\Q^n}$ is ambiguous in each of these cases and the `gap' is actually quite wide. In the most famous case of the plane we have $\chi(T_{\Q^2})=2$ while $\chi_m(T_{\Q^2})\geq5 $.

It turns out that in general for higher dimensions Proposition \ref{usecor} provides a stronger bound on $\chi_m(T_{\K^n})$ than Theorem \ref{main}. For $n=5,\dots,12$ the known bounds on $\chi(\R^n)$ are better than $n+3$ \cite{kuparaigo}. What's more, it is known that $\chi(\R^n)$ grows exponentially with $n$ \cite{raigo}, so Proposition \ref{usecor} will be stronger than Theorem \ref{main} for all subsequent $n$. However, for $n\geq5$ we know of no colourings of $\Q^n$ (or $\K^n$) which provide further ambiguous examples. Raigorodskii's survey \cite[p.111]{raigo} suggests that it is known that $\chi(\Q^5)\leq8$, citing Chilakamarri \cite{chilaka}. However, Chilakamarri only conjectures that $\chi(\Q^5)=8$, and we were unable to find any proof of this proposition elsewhere in the literature. Interestingly, Cibulka has recently shown that $\chi(\Q^5)\geq 8$ \cite{cibulka}, so along with Cantwell's result that $\chi(\R^5)\geq 9$ \cite{cantwell}, an 8-colouring of $\Q^5$ would prove Chilakamarri's conjecture and furnish a further ambiguous example in $T_{\Q^5}$.

Concentrating now on dimension $2$, there are some other fields $\K$ for which useful results about $\chi(\K^2)$ are known. The following results concern quadratic extensions $\Q[\sqrt{n}]$ where $n$ is a positive square free integer. Johnson \cite{johnson} showed that $\chi(\Q[\sqrt{n}]^2)=2$ for $n \equiv_4 1,2$ and Fischer \cite{fischer} showed that $\chi(\Q[\sqrt{n}]^2)\leq3$ for $n \equiv_3 0,1$ and that $\chi(\Q[\sqrt{n}]^2)\leq4$ for $n \equiv_8 3$, so for all these cases $T_{\K^2}$ has ambiguous chromatic number. In particular we have the example $T_{\Q[\sqrt{3}]^2}$ which contains all equilateral triangles which have some edge vector in $\Q^2$, and hence many copies of the triangle lattice. In this case we have $\chi = 3$ and $\chi_m \geq 5$.

\subsection*{Acknowledgements}
Thanks to Lashi Bandara, Burkard Polster, Marty Ross, Moritz Schmitt, Ian Wanless and Günter Ziegler for their helpful advice, and especially to Boris Bukh for his suggestions and criticism. Thanks also to the referee for pointing to some relevant recent references.

\bibliographystyle{siam}

\bibliography{references}

\def\cprime{$'$}
\begin{thebibliography}{10}

\bibitem{benper}
{\sc M.~Benda and M.~Perles}, {\em Colorings of metric spaces}, Geombinatorics,
  9 (2000), pp.~113--126.

\bibitem{cantwell}
{\sc K.~Cantwell}, {\em Finite {E}uclidean {R}amsey theory}, J. Combin. Theory
  Ser. A, 73 (1996), pp.~273--285.

\bibitem{chilaka}
{\sc K.~B. Chilakamarri}, {\em The unit-distance graph problem: a brief survey
  and some new results}, Bull. Inst. Combin. Appl., 8 (1993), pp.~39--60.

\bibitem{cibulka}
{\sc J.~Cibulka}, {\em On the chromatic number of real and rational spaces},
  Geombinatorics, 18 (2008), pp.~53--65.

\bibitem{croftlat}
{\sc H.~T. Croft}, {\em Three lattice-point problems of {S}teinhaus}, Quart. J.
  Math. Oxford Ser. (2), 33 (1982), pp.~71--83.

\bibitem{EdB}
{\sc N.~G. de~Bruijn and P.~Erd\H{o}s}, {\em A colour problem for infinite
  graphs and a problem in the theory of relations}, Nederl. Akad. Wetensch.
  Proc. Ser. A. {\bf 54} = Indagationes Math., 13 (1951), pp.~369--373.

\bibitem{falconer}
{\sc K.~J. Falconer}, {\em The realization of distances in measurable subsets
  covering {${\mathbb R}\sp{n}$}}, J. Combin. Theory Ser. A, 31 (1981),
  pp.~184--189.

\bibitem{fischer}
{\sc K.~G. Fischer}, {\em Additive {$K$}-colorable extensions of the rational
  plane}, Discrete Math., 82 (1990), pp.~181--195.

\bibitem{johnson}
{\sc P.~D. Johnson, Jr.}, {\em Two-colorings of real quadratic extensions of
  {$Q\sp 2$} that forbid many distances}, Congr. Numer., 60 (1987), pp.~51--58.
\newblock Eighteenth Southeastern International Conference on Combinatorics,
  Graph Theory, and Computing (Boca Raton, Fla., 1987).

\bibitem{kuparaigo}
{\sc A.~Kupavskii and A.~Raigorodskii}, {\em On the chromatic numbers of
  small-dimensional euclidean spaces}, Electronic Notes in Discrete
  Mathematics, 34 (2009), pp.~435 -- 439.
\newblock European Conference on Combinatorics, Graph Theory and Applications
  (EuroComb 2009).

\bibitem{raigo}
{\sc A.~M. Ra{\u\i}gorodski{\u\i}}, {\em Borsuk's problem and the chromatic
  numbers of some metric spaces}, Russian Math. Surveys, 56 (2001),
  pp.~103--139.

\bibitem{schmutz}
{\sc E.~Schmutz}, {\em Rational points on the unit sphere}, Cent. Eur. J.
  Math., 6 (2008), pp.~482--487.

\bibitem{ss03}
{\sc S.~Shelah and A.~Soifer}, {\em Axiom of choice and chromatic number of the
  plane}, J. Combin. Theory Ser. A, 103 (2003), pp.~387--391.

\bibitem{soifer05}
{\sc A.~Soifer}, {\em Axiom of choice and chromatic number of {$\R^n$}}, J.
  Combin. Theory Ser. A, 110 (2005), pp.~169--173.

\bibitem{soifbook}
\leavevmode\vrule height 2pt depth -1.6pt width 23pt, {\em The mathematical
  coloring book}, Springer, New York, 2009.
\newblock Mathematics of coloring and the colorful life of its creators, With
  forewords by Branko Gr{\"u}nbaum, Peter D. Johnson, Jr. and Cecil Rousseau.

\bibitem{ss04}
{\sc A.~Soifer and S.~Shelah}, {\em Axiom of choice and chromatic number:
  examples on the plane}, J. Combin. Theory Ser. A, 105 (2004), pp.~359--364.

\bibitem{sakey}
{\sc L.~A. Sz{\'e}kely}, {\em Measurable chromatic number of geometric graphs
  and sets without some distances in {E}uclidean space}, Combinatorica, 4
  (1984), pp.~213--218.

\bibitem{woodall}
{\sc D.~R. Woodall}, {\em Distances realized by sets covering the plane}, J.
  Combin. Theory Ser. A, 14 (1973), pp.~187--200.

\end{thebibliography}

\end{document}